\documentclass[journal]{IEEEtran}

\usepackage{cite}
\usepackage{amsmath,amssymb,amsfonts,bm}
\usepackage{bbm}
\usepackage{graphicx}
\usepackage{subfigure}
\usepackage{textcomp}
\usepackage{algorithm}
\usepackage{algpseudocode}

\usepackage{xcolor}
\usepackage{url}

\newtheorem{theorem}{Theorem}
\newtheorem{lemma}{Lemma}
\newtheorem{assum}{Assumption}

\newtheorem{remark}{Remark}

\newenvironment{proof}{\begin{IEEEproof}}{\end{IEEEproof}}

\begin{document}
	
	\title{Convergence and Sample Complexity of Policy Gradient Methods for Stabilizing Linear Systems
		\thanks{The work was supported by National Natural Science Foundation of China under Grant no. 62033006 and Tsinghua University Initiative Scientific Research Program. ({\em Corresponding author: Keyou You})}
		\thanks{The authors are with the Department of Automation and BNRist, Tsinghua University, Beijing 100084, China. e-mail: zhaofr18@mails.tsinghua.edu.cn, fxy20@mails.tsinghua.edu.cn, youky@tsinghua.edu.cn. }
		\author{Feiran Zhao, Xingyun Fu, Keyou You, \IEEEmembership{Senior Member,~IEEE}}}
	\maketitle
	
	\begin{abstract} System stabilization via policy gradient (PG) methods has drawn increasing attention in both control and machine learning communities. In this paper, we study their convergence and sample complexity for stabilizing linear time-invariant systems in terms of the number of system rollouts. Our analysis is built upon a discounted linear quadratic regulator (LQR) method which alternatively updates the policy and the discount factor of the LQR problem. Firstly, we propose an explicit rule to adaptively adjust the discount factor by exploring the stability margin of a linear control policy. Then, we establish the sample complexity of PG methods for stabilization, which only adds a coefficient logarithmic in the spectral radius of the state matrix to that for solving the LQR problem with a prior stabilizing policy. Finally, we perform simulations to validate our theoretical findings and demonstrate the effectiveness of our method on a class of nonlinear systems.
	\end{abstract}
	
\section{Introduction}
Recent years have witnessed tremendous success in applications of reinforcement learning (RL) in continuous control \cite{tobin2017domain, andrychowicz2020learning, recht2019tour} and sequential decision-making problems~\cite{mnih2015human-level, silver2016mastering}. Most of these successes hinge on the developments of the policy gradient (PG) method, an essential approach in modern RL, which directly searches over the policy space to minimize a cost function. Such an end-to-end approach only requires samples of the cost from system rollouts, and thus is conceptually simple to implement. In fact, the PG method has been directly applied to a wide variety of tasks for complex dynamics~\cite{levine2016end}.

The success of RL has inspired an increasing interest in studying the convergence and sample complexity of PG methods on classical control problems, which serve as ideal benchmarks to understand their efficiency and complexities  \cite{fazel2018global, bu2019lqr, malik2019derivative, mohammadi2020linear, mohammadi2021convergence, zhang2021policy,zhao2021global, zhang2019policy, zheng2021analysis, duan2022optimization,furieri2020learning, li2021distributed, fatkhullin2021optimizing,  gravell2020learn, tu2019gap, yang2019provably}. The seminal work of \cite{fazel2018global} first shows that PG methods have global convergence guarantees for the celebrated linear quadratic regulator (LQR) problem. Then, sample complexities of PG methods in terms of the number of system rollouts are established for both discrete-time \cite{malik2019derivative, mohammadi2020linear} and continuous-time LQR \cite{mohammadi2021convergence}, and they have also been applied to many variants of LQR problems, such as risk-sensitive/constrained LQ control \cite{zhang2021policy,zhao2021global}, LQ games \cite{zhang2019policy}, linear quadratic Gaussian (LQG) \cite{zheng2021analysis, duan2022optimization} and decentralized control \cite{furieri2020learning, li2021distributed}. Though these advances lead to fruitful and profound results for PG methods, they all require a common assumption: a stabilizing policy must be known \textit{a prior}. This  is essential to PG methods~\cite{fazel2018global, zhang2021policy,zhang2019policy,zheng2021analysis}, since their iterative local search makes sense \textit{only if the cost is finite}. In fact, learning a stabilizing policy from samples of the cost is a non-trivial task, and even of the same importance as solving the LQR problem \cite{hu2022sample, kang2023minimum}. 

Very recently, \textit{discount methods} have been introduced to address the stabilization problem for linear systems \cite{lamperski2020computing, jing2021learning} and nonlinear systems~\cite{perdomo2021stabilizing}. The discount method is referred to as a class of iterative methods that solve a sequence of discounted LQR problems with increasing discount factors, and is originally proposed for escaping local optimal policies in multi-agent control systems~\cite{feng2020escaping,feng2021damping}. However, Ref. \cite{lamperski2020computing, jing2021learning} do not provide finite-time convergence guarantees to learn stabilizing policies. While Ref. \cite{perdomo2021stabilizing} achieves stabilization in a finite number of iterations on the discount factors, it does not show the sample complexity even for linear systems.


The sample complexity in \cite{perdomo2021stabilizing} cannot be explicitly provided since (a) their discount method requires to completely solve the discounted LQR problem at each iteration, and (b) the update of the discount factor relies on a binary or random search, which involves an unpredictable number of iterations. In contrast, we propose a new discount method where (a) we only need to reduce the cost to a uniform level, and (b) we design an explicit update rule to recursively adjust the discount factor. Our key is to establish a connection between the stability margin of a linear control policy and its discounted LQR cost by using Lyapunov theory, which can be viewed as a natural generalization of the stability margin of the LQR~\cite{safonov1977gain}. Then, we provide performance and finite-time convergence guarantees for our discount method. 
To show its sample complexity, we adopt zero-order PG methods to reduce the cost and estimate the update rate of the discount factor via a simulator used to return samples of the cost~\cite{bertsekas2019reinforcement}. Finally, we validate our theoretical results via simulations and verify that our methods can also efficiently stabilize nonlinear dynamical systems around the equilibrium point. 

Our results are the first to reveal that the sample complexity of PG methods for stabilizing linear systems only adds a coefficient logarithmic in the spectral radius of the state matrix to that for solving the LQR problem with a prior stabilizing policy~\cite{mohammadi2021convergence}. Moreover, our explicit update rule can be easily used to study sample complexity for other RL based discount methods, e.g.,~\cite{jing2021learning,lamperski2020computing}.

To clarify our contributions in the literature, we further compare with recent methods for system stabilization, which can be broadly categorized as model-based approach, direct data-driven control, and RL approach.
 
\textbf{Model-based approach.} It learns a descriptive model via system identification using samples from  multiple trajectories~\cite{dean2020sample, zheng2020non} or a single trajectory ~\cite{abbasi2011regret, ibrahimi2012efficient, talebi2021regularizability, chen2021black, lale2020explore, faradonbeh2018finite, hu2022sample}. In this context, the sample complexity is quantified in terms of the number of state-input pairs for system identification (SysID). A necessary condition for SysID is the persistent excitation (PE) condition, which requires at least $n+m$ samples with $n$ and $m$ being the state and input dimension, respectively. Consequently, their sample complexity typically scales linearly in $n$. While a recent work \cite{hu2022sample} improves the complexity to be linear in the number of unstable eigenvalues of the state matrix, they only focus on linear systems with diagonalizable state matrix. In this work, our goal is to strengthen the theoretical foundation of PG methods. Since we only have access to a scalar cost per system rollout, our sample complexity is quadratic in $n$ in terms of the number of system rollouts. Thus, it is comparable with model-based approach where each sample contains a state vector with dimension $n$.

\textbf{Direct data-driven control.} There is another line of research that directly synthesizes controllers from a given set of system trajectories. The seminal work \cite{willems2005note} establishes Willems' fundamental lemma  to show that all the trajectories of a linear system can be represented by a linear combination of historical trajectories under the PE condition. Then, a data-based semi-definite program (SDP) is designed in \cite{de2019formulas, de2021low} to find a stabilizing controller when the additive noise of the system dynamics satisfies a quadratic matrix inequality. Without the PE condition, Ref. \cite{van2020data} provides a  necessary and sufficient condition for stabilizing noiseless linear systems with state feedback. For the stochastic case, a so-called matrix $\mathcal{S}$-lemma is proposed in \cite{van2020noisy, bisoffi2021trade} to find a common stabilizing controller of a set of admissible systems consistent with data samples. While this approach provides a novel perspective on the stabilization problem, it relies on the solution to an SDP, the existence of which is sensitive to the data quality~\cite{van2020noisy}. Different from this work, the above references cannot provide any sample complexity in finding a stabilizing controller. In fact, a recent work~\cite{kang2023minimum} shows that the sample complexity of stabilization via direct data-driven approaches is almost identical to that of SysID, which is $n+m$ in terms of the number of state-input pairs. We refer the readers to \cite{markovsky2021behavioral} for a comprehensive review of this direction. 

\textbf{RL approach.} Instead of identifying an explicit dynamical model or solving an SDP, the RL approach directly updates the policy using samples~\cite{perdomo2021stabilizing, lamperski2020computing, chen2022homotopic}. Ref. \cite{lamperski2020computing} is the first to combine the discount method and Q-learning to stabilize deterministic linear systems, while they require to estimate the Q-function using the states and control inputs. In contrast, our PG method does not need to identify any explicit  Q-function and only requires the cost of a system rollout from the simulator. Moreover, using our explicit update rule for the discount factor, we can provide finite-time convergence guarantees for their method. A recent work \cite{chen2022homotopic} considers the continuous-time linear systems and learns a stabilizing controller by a so-called homotopic policy iteration scheme. This approach gradually pushes the unstable poles to the left-half complex plane with a lower bounded distance each iteration, instead of multiplying a damping factor. Again, they fail to prove finite-time convergence guarantees as well as the sample complexity.

The rest of this paper is structured as follows. In Section \ref{sec:problem}, we describe the stabilization problem of linear systems with randomized initial states. In Section \ref{sec:learning}, we propose our explicit discount method algorithm where the update rule for the discount factor is designed by studying the stability margin of a linear control policy. Section \ref{sec:sample} applies PG methods for our discount method and shows the sample complexity in terms of the number of rollouts. Section \ref{sec:exp} performs simulations on both linear and nonlinear systems to verify the effectiveness of the proposed method. In Section \ref{sec:conclusion}, we draw some concluding remarks and discuss future directions.

\textbf{Notations.} We use $\rho(\cdot)$ to denote the spectral radius of a matrix. We use $\|\cdot\|_{\text{F}}$ and $\|\cdot\|$ to denote the Frobenius norm and $2$-norm of a matrix, respectively. Let $\underline{\sigma}(\cdot)$ be the minimal singular value of a matrix. $\text{Tr}(\cdot)$ denotes the trace function. Let $I_{n\times n}$ be the $n$-by-$n$ identity matrix. Let $\mathcal{U}^{d-1} \subset \mathbb{R}^d$ be the unit sphere of dimension $d-1$. We use $\mathcal{O}(\epsilon)$ to denote some constant proportional to $\epsilon$. We use $\widetilde{\mathcal{O}}(\cdot)$ to neglect the logarithmic term in the dimension $n$ of $\mathcal{O}(\cdot)$.

\section{Problem formulation}\label{sec:problem}

Consider the following discrete-time linear time-invariant system
\begin{equation}\label{equ:sys}
x_{t+1} = Ax_t + Bu_t, ~~x_0 \sim \mathcal{D},
\end{equation}
where $x_t \in \mathbb{R}^n$ is the state vector and $u_t\in \mathbb{R}^m$ is the control input vector. The system $(A,B)$ is assumed to be controllable as in \cite{fazel2018global,bu2019lqr,malik2019derivative}. The initial state $x_0$ is sampled from a distribution $\mathcal{D}$, on which we make the following mild assumption.

\begin{assum}
	\label{assumption:D}
	The distribution $\mathcal{D}$ has zero mean and covariance $\mathbb{E}[x_0x_0^{\top}] = I$ with bounded support, i.e., $\|x_0\|\leq d$ for some constant $d>0$. 
\end{assum}

{Assumption \ref{assumption:D} is made only for simplicity of theoretical analysis, and our results can be extended to the unbounded distributions with strong tail decay, e.g., Gaussian distributions, by using the high-probability bounds and standard truncation arguments in [9]. In fact, we have validated our main results over the Gaussian distribution in the simulation.
}

The policy gradient (PG) method parameterizes the control policy with a feedback gain matrix $K$, and uses gradient descent method to solve the following discounted LQR problem
\begin{equation}\label{prob:dlqr}
\begin{aligned}
&\mathop{\text{ minimize}}\limits_K ~ J_{\gamma}(K) := \mathbb{E}_{x_0} \sum_{t=0}^{\infty}\gamma^{t} (x_{t}^{\top} Q x_{t}+u_{t}^{\top} R u_{t})~~~~\\
&~\text {subject to} ~(\ref{equ:sys})~\text{and}~u_t = -K x_t,
\end{aligned}
\end{equation}
where $0<\gamma \leq 1$ is the discount factor, and $Q>0, R>0$ are user-specified penalty matrices. {Here, the positive definiteness of $Q$ is made only for technical simplicity, and can be relaxed to $Q\geq 0$ using the techniques in \cite{bu2020global}.} The PG is computed using $(A,B)$ in the model-based setting or estimated from system rollouts in the sample-based setting (c.f. Section \ref{sec:sample}). As an iterative search method, an initial policy $K^0$ is required to render finite cost, i.e., $J_{\gamma}(K^0)<\infty$.


While PG methods have successfully shown theoretical guarantees in many LQ control problems~\cite{mohammadi2021convergence, zhang2021policy, zhao2021global,zhang2019policy}, their common and essential assumption is the prior knowledge of a stabilizing gain $K$, i.e., $\rho(A-BK)<1$. This is non-trivial to achieve for unknown $(A,B)$ and how to remove it has been acknowledged as an open problem in \cite{fazel2018global}. To address it, a discount method has been proposed in \cite{perdomo2021stabilizing}. Yet it is unable to quantify its sample complexity, rendering the efficiency of PG methods still unclear.

{In this paper, we study the convergence and sample complexity of PG methods for stabilizing the linear system in \eqref{equ:sys} under the sample-based setting. In particular, we only have access some prior knowledge of $\rho(A)$ and a simulator which returns samples of the cost of \eqref{prob:dlqr}.  Note that such a simulator has been widely used  in the literature of PG methods, e.g., \cite{fazel2018global, bu2019lqr, malik2019derivative, mohammadi2020linear, mohammadi2021convergence, perdomo2021stabilizing}.}
Particularly, we propose a new discount method with an explicit update rule to adaptively adjust the discount factor. Then, we adopt the PG method to implement it via a simulator and rigorously quantify its  sample complexity.

\section{Stabilization via Discount methods}\label{sec:learning}

In this section, we first provide some preliminary results on the discounted LQR problem and recapitulate the discount method. Then, we design an explicit update rule for the discount factor by studying the stability margin of linear control policies, based on which we are able to achieve convergence guarantees for finding a stabilizing policy.

\subsection{The discount method}
We first provide some important properties of the discounted LQR problem (\ref{prob:dlqr}). It is well-known that the discounted cost in (\ref{prob:dlqr}) is equivalent to the standard LQR cost of the following damped dynamical system \cite{malik2019derivative,perdomo2021stabilizing}
\begin{equation}\label{equ:rescale_dyna}
x_{t+1} = \sqrt{\gamma}(A - BK) x_t,~x_0 \sim \mathcal{D},
\end{equation}
i.e., it holds that
\begin{equation}\label{equ:equivalent}
\begin{aligned}
&J_{\gamma}(K) = \mathbb{E}_{x_0} \sum_{t=0}^{\infty} (x_{t}^{\top} Q x_{t}+u_{t}^{\top} R u_{t})\\
&~\text {subject to} ~(\ref{equ:rescale_dyna})~\text{and}~u_t = -K x_t.
\end{aligned}
\end{equation}

We now provide a condition for the cost $J_{\gamma}(K)$ to be finite. To this end, define the stability region of the damped system (\ref{equ:rescale_dyna}) as
$$
\mathcal{S}_{\gamma} = \{K | \sqrt{\gamma}\rho(A-BK) < 1 \}.
$$ 
\begin{lemma}[\cite{fazel2018global}]\label{lem:closed} Consider the discounted LQR problem (\ref{prob:dlqr}). Then,
$J_{\gamma}(K) < \infty$ if and only if $K \in \mathcal{S}_{\gamma}$. Moreover, the finite cost has a closed-form expression 
	$$J_{\gamma}(K) = \text{Tr}(P(\gamma))=\text{Tr}((Q+K^{\top}RK)\Sigma_K(\gamma)),$$
	where $P(\gamma)$ and $\Sigma_K(\gamma)$ are positive semi-definite solutions to the Lyapunov equations
	\begin{align}
	&P(\gamma) = Q + K^{\top}RK + \gamma (A-BK)^{\top}P(\gamma)(A-BK), \label{equ:Pi} \\
	&\Sigma_K(\gamma) = I + \gamma (A-BK)\Sigma_K(\gamma)(A-BK)^{\top}, \label{def:cov}
	\end{align}
	respectively. 
\end{lemma}

For simplicity of notation, we use the shorthand $\Sigma:= \Sigma_K(\gamma)$ when $K$ and $\gamma$ are clear from the context.

\begin{figure}[t]
	\centering
	\includegraphics[height=46mm]{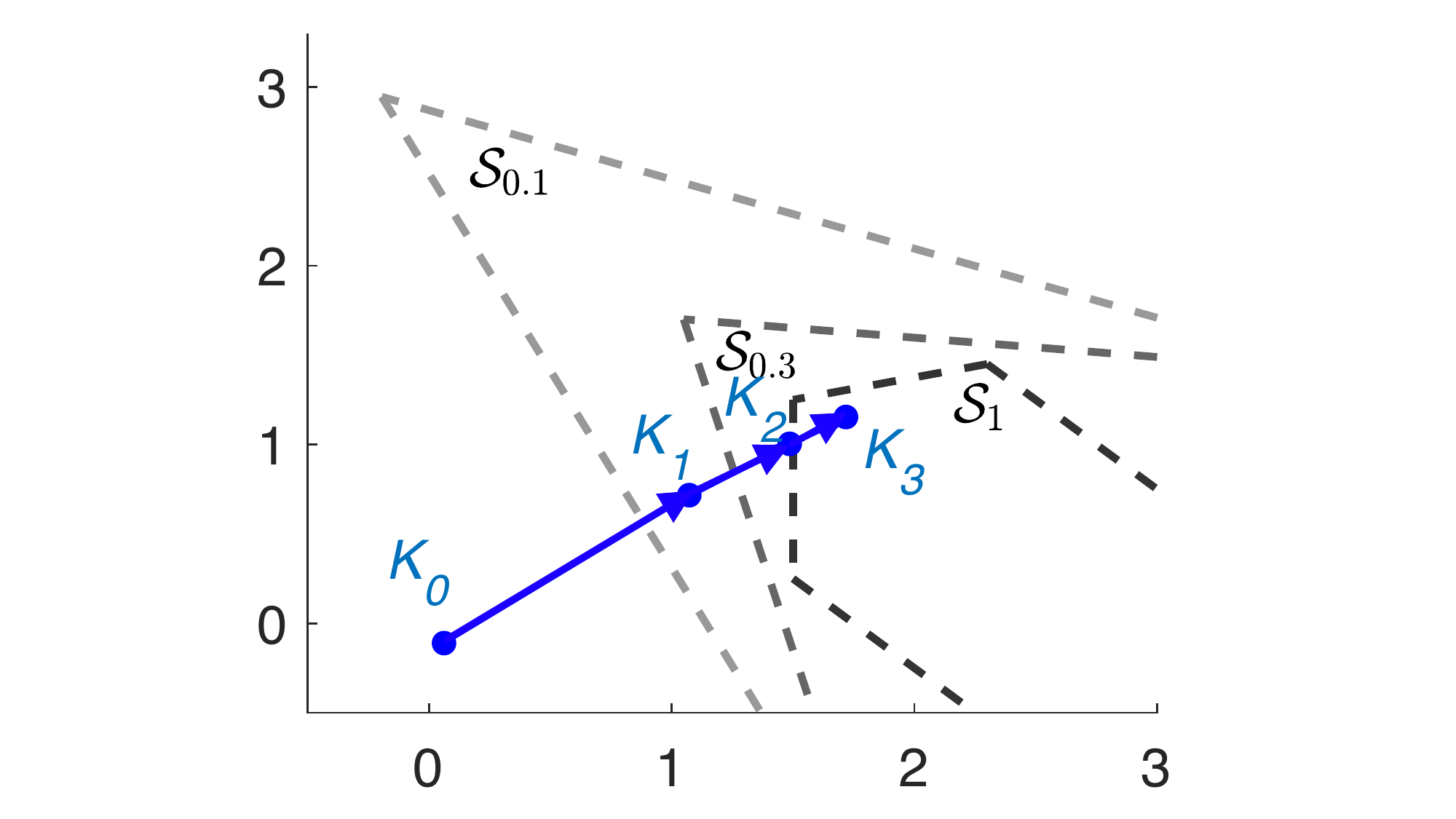}
	\caption{Policy trajectory of the proposed method in a 2-dimensional example provided in Section \ref{sec:exp}.  We plot the boundaries of the sets $\mathcal{S}_{\gamma} = \{K | \sqrt{\gamma}\rho(A-BK) < 1 \}$ with $\gamma = 0.1, 0.3, 1$, and select three intermediate policies denoted by $K_1,K_2,K_3$ for illustration. From an initial policy $K_0$, our algorithm drives the policy towards the stabilizing set $\mathcal{S}_{1}$.	
	}
	\label{pic:pgland}
\end{figure}

The discount methods in \cite{lamperski2020computing, perdomo2021stabilizing} alternate between the following two procedures
\begin{subequations}\label{equ:disc}
	\begin{align}
	&\text{Solve}~K^k \in \arg \min_K J_{\gamma^k}(K), \label{equ:disc1}\\
	&\text{Search}~\gamma^{k+1}>\gamma^k \text{such that} \sqrt{\gamma^{k+1}}\rho(A-BK^k)<1.\label{equ:disc2}
\end{align}
\end{subequations}
Given a fixed $\gamma^k$, the step \eqref{equ:disc1} finds an optimal gain $K^k$ of the discounted LQR problem. By the well-known stability margin of the LQR~\cite{safonov1977gain}, the optimal gain $K^k$ is bounded away from the boundary of $\mathcal{S}_{\gamma^k}$, i.e., there exists $\gamma^{k+1}>\gamma^{k}$ such that $\sqrt{\gamma^{k+1}}\rho(A-BK^k)<1$. Then, a binary or random search method is used to find such a discount factor $\gamma^{k+1}$ in (\ref{equ:disc2}). For example, a binary search is proposed in \cite{perdomo2021stabilizing} to find $\gamma^{k+1} \in (\gamma^k, 1]$ such that
$$
2.5 J_{\gamma^k}(K^{k+1}) \leq J_{\gamma^{k+1}}(K^{k+1}) \leq 8J_{\gamma^k}(K^{k+1}).
$$
However, \eqref{equ:disc1} requires to minimize the discounted LQR cost per iteration, which is time-consuming and challenging to implement in practice. Moreover, the binary or random search in \eqref{equ:disc2} may result in an unpredictable number of iterations.

Since their discount factor in (\ref{equ:disc2}) is implicitly computed, we refer to (\ref{equ:disc}) as the \textit{implicit} discount method. In this paper, we propose an \textit{explicit} discount method, alternating as
\begin{subequations}\label{equ:expmethd}
	\begin{align}
	&\text{Find}~K^k ~\text{such that}~J_{\gamma^k}(K^k)<D, \label{equ:exp1}\\
	&\text{Update}~\gamma^{k+1}= (1+\xi\alpha^k)\gamma^k,\label{equ:exp2}
	\end{align}
\end{subequations}
where $D, 0<\xi<1$ are two constants to be designed later, and $\alpha^k>0$ is the {update rate} of the discount factor that can be explicitly computed using the value of $J_{\gamma^k}(K^k)$. Instead of minimizing the cost in (\ref{equ:disc1}), we only need to reduce it to a uniform level $D$ in (\ref{equ:exp1}). Moreover, (\ref{equ:exp2}) is an explicit update rule for the discount factor. Fig. \ref{pic:pgland} provides an illustration of our method.

\subsection{The explicit update rule for the discount factor}
To specify the update rule in (\ref{equ:exp2}), we exploit the relation between the stability margin and the cost of a linear control policy. Different from the classical stability margin of the {LQR} \cite{safonov1977gain}, we study the following problem regarding the stability margin of an \textit{arbitrary} policy $K$:

\textit{Given a policy $K\in \mathcal{S}_{\gamma}$, explicitly find a larger discount factor $\gamma'>\gamma$ such that $K\in \mathcal{S}_{\gamma'}$ without the knowledge of $(A,B)$. }

We first provide a sufficient condition on $\gamma'$ using the Lyapunov equation.
\begin{lemma}\label{theo:rule}
	If $K \in \mathcal{S}_{\gamma}$ and $\gamma'$ satisfies
	\begin{equation}\label{equ:rule}
	(1-\gamma/\gamma')P < Q + K^{\top}RK,
	\end{equation}
	where $P$ is a positive semi-definite solution of the Lyapunov equation in (\ref{equ:Pi}), then $K \in \mathcal{S}_{\gamma'}$.
\end{lemma}

\begin{proof}
	Consider the following damped system
	\begin{equation}\label{sys:lyapnov}
	x_{t+1} = \sqrt{\gamma'}(A-BK)x_t.
	\end{equation}	
	Define $V(x) = x^{\top}Px$. Clearly, $V(x)$ satisfies that $V(0) = 0$ and $V(x)>0,\forall x\neq 0$. Moreover, we obtain 
	\begin{align*}
		&V(x_{t+1}) - V(x_t) \\
		&= \gamma'x_t^{\top}(A-BK)^{\top}P(A-BK)^{\top}x_t - x_t^{\top}P x_t \\
		&= x_t^{\top}\left(\frac{\gamma'}{\gamma}(P-Q-K^{\top}RK) - P\right) x_t,
	\end{align*}
	where the last equality follows from (\ref{equ:Pi}). 
	
	Jointly with (\ref{equ:rule}), it follows that $V(x_{t+1}) - V(x_t)<0$, which verifies that $V(x) = x^{\top}Px$ is a Lyapunov function for (\ref{sys:lyapnov}). Thus, the damped system \eqref{sys:lyapnov} is asymptotically stable. 
\end{proof}

However, it is impossible to directly obtain $\gamma'$ by Lemma \ref{theo:rule} as the computation of $P$ involves the unknown $(A,B)$. Next, we use Lemma \ref{lem:closed}  to show that the condition (\ref{equ:rule}) can be achieved using the cost $J_{\gamma}(K)$.

\begin{theorem}\label{coro}
	If $K \in \mathcal{S}_{\gamma}$ and $\gamma' \leq (1+\alpha)\gamma$ with
	\begin{equation}\label{equ:coro}
	\alpha = \frac{\underline{\sigma}(Q+K^{\top}RK)}{J_{\gamma}(K) - \underline{\sigma}(Q+K^{\top}RK)},
	\end{equation}
	then $K \in \mathcal{S}_{\gamma'}$.
\end{theorem}

\begin{proof}
	To satisfy \eqref{equ:rule}, it suffices for $\gamma'$ to satisfy
	$$
	1-{\gamma}/{\gamma'} < {\underline{\sigma}(Q+K^{\top}RK)}/{\text{Tr}(P) }.
	$$
	The proof follows from Lemma \ref{lem:closed} that
	$
	J_{\gamma}(K) =\text{Tr}(P).
	$
\end{proof}

Theorem \ref{coro} reveals how much the system matrices $A$ and $B$ can be scaled to maintain stability of $\sqrt{\gamma}(A-BK)$ for a given policy $K$. Since $\alpha \geq \underline{\sigma}(Q)/(J_{\gamma}(K)-\underline{\sigma}(Q))$ in (\ref{equ:coro}), the scaling is proportional to the inverse of the cost $J_{\gamma}(K)$. Thus, we only need to reduce the cost to $D$ in (\ref{equ:exp1}) to ensure the uniform lower bound of $\alpha^k$. Moreover, Theorem \ref{coro} suggests an explicit update rule for the discount factor in the form of (\ref{equ:exp2}), the computation of which is trivial since $J_{\gamma}(K)$ can be evaluated from a simulator and $\underline{\sigma}(Q+K^{\top}RK)$ is known.

\begin{remark}\label{remark1}
	Herein we only consider that $A$ and $B$ are scaled, which is sufficient for our discount method, and use the discount factor to quantify the stability margin. How to rigorously evaluate its connection to the classical notion of stability margin is beyond the scope of this paper.\end{remark}

In the sequel, we propose the explicit discount method in the form of \eqref{equ:expmethd} to find a stabilizing policy.

\begin{algorithm}[t]
	\caption{The explicit discount method}
	\label{alg:exact}
	\begin{algorithmic}[1]
		\Require
		An initial policy $K^0 = 0$, a discount factor $\gamma^0<1/\rho^2(A)$, a bound $D$, and a coefficient $0<\xi<1$.
		\For{$k=0,1,\cdots$}
		\State \textbf{Cost reduction:} find $K^{k+1}$ such that
		\begin{equation}\label{equ:cost_req}
		J_{\gamma^k}(K^{k+1}) \leq D.
		\end{equation}
		\State \textbf{Update of the discount factor:}  
		\begin{equation}\label{equ:rule2}
		\begin{aligned}
		&\gamma^{k+1} = (1+\xi\alpha^k)\gamma^k ~\text{with} \\
		&\alpha^k = \frac{\underline{\sigma}(Q+(K^{k+1})^{\top}RK^{k+1})}{J_{\gamma^k}(K^{k+1}) - \underline{\sigma}(Q+(K^{k+1})^{\top}RK^{k+1})}.
		\end{aligned}
		\end{equation} 
		\If{$\gamma^{k+1} \geq  1 $} return a stabilizing policy $K^{k+1}$.
		\EndIf
		\EndFor
	\end{algorithmic}
\end{algorithm}

\subsection{The explicit discount method with finite-time convergence guarantees}
Our discount method in the form of \eqref{equ:expmethd} is presented in Algorithm \ref{alg:exact}. Since PG methods require an initial policy with finite cost, we initialize $K^0 = 0$ and $\gamma^0 \leq 1/\rho^2(A)$ such that $K^0 \in \mathcal{S}_{\gamma^0}$ under the prior knowledge of $\rho(A)$.

In the cost reduction step, we use PG methods to reduce the discounted LQR cost to a uniform level $D$ starting from an initial policy $K^k \in \mathcal{S}_{\gamma^k}$, the details of which will be specified in the next section. This step ensures a uniform lower bound of the update rate of the discount factor. The selection of $D$ is important to achieve (\ref{equ:cost_req}). Particularly, $D$ must be chosen such that $D>J^*_{\gamma^k}$ at the $k$-th iteration, where $J^*_{\gamma^k}$ is the optimal cost by definition $J^*_{\gamma}:= \min_K J_{\gamma}(K).$ Thanks to the following lemma, it suffices to set $D > J^*_1$. 
\begin{lemma}\label{lem:compare}
	For $0<\gamma_1<\gamma_2\leq 1$, it holds that 
	$
	J^*_{\gamma_1} < J^*_{\gamma_2}.
	$
\end{lemma}

While the design of $D$ needs some knowledge on the upper bound of $J^*_1$, we can set $D$ to be sufficiently large in practice to circumvent this requirement with the sacrifice of a longer running time, as to be shown in Theorem \ref{thm:tech}. The cost requirement in (\ref{equ:cost_req}) is more tractable than the criteria $J_{\gamma^k}(K^{k+1}) - J_{\gamma^k}^*\leq \epsilon$ in \cite{perdomo2021stabilizing}, which is difficult to evaluate since $J_{\gamma^k}^*$ is unknown.

The update rule of the discount factor is provided in (\ref{equ:rule2}), where we additionally multiply a coefficient $\xi\in(0,1)$ to the update rate $\alpha^k$ to 
provide performance guarantees for the resulted policy. Particularly, it ensures that $\sqrt{\gamma^{k+1}}\rho(A-BK^{k+1})$ is ``strongly" stable for $k \in \mathbb{N}$, i.e., the spectral radius of the damped closed-loop matrix is bounded away from one. 
By the cost reduction step, the update rate has a lower bound depending on $D$, leading to the following finite-time convergence guarantees of Algorithm \ref{alg:exact}.
\begin{theorem}\label{thm:tech}
	If $D > J^*_1$ and $$k\ge {(D-\underline{\sigma}(Q))}\log ({1}/{\gamma^0})/(\xi\underline{\sigma}(Q)),$$ then Algorithm \ref{alg:exact} returns a stabilizing policy $K^k$ with $\rho(A-BK^k)<\sqrt{1-(1-\xi)\underline{\sigma}(Q)/D}$.
\end{theorem}

The proof of Theorem \ref{thm:tech} is provided in Appendix \ref{app:Lemma 5}. Theorem \ref{thm:tech} explicitly shows the dependence of the number of iterations and performance of Algorithm \ref{alg:exact} on the user-specified parameters $\gamma^0, D$, and $\xi$. Particularly, the upper bound of $\rho(A-BK^k)$ significantly relies on the coefficient $\xi$. If $\xi = 1$, then the resulted closed-loop matrix is only asymptotically stable and it is impossible to find a uniform bound for $\rho(A-BK^k)$. 
Theorem \ref{thm:tech} also shows a salient improvement over \cite{perdomo2021stabilizing} in the number of total iterations, i.e., their algorithm needs $64 (J_{1}^*)^4 \log ({1}/{\gamma^0})$ iterations to converge, which is a $4$-th order polynomial of ours.


In the next section, we use a simulator to apply PG methods with gradient estimate to perform the cost reduction step of Algorithm \ref{alg:exact} and estimate the update rate in \eqref{equ:rule2} as $J_{\gamma^k}(K^{k+1})$ is unknown.

\section{Sample complexity of PG based discount methods for stabilization}\label{sec:sample}
In this section, we use a simulator to return samples of the cost in (\ref{equ:equivalent}) given $\gamma$ and $K$ from system rollouts. Particularly, we show how to use the simulator to estimate the PG for cost reduction and the update rate of the discount factor in (\ref{equ:rule2}). Then, we show the sample complexity of our PG based discount method in terms of the number of system rollouts to learn a stabilizing policy, where a system rollout returns a sample of the cost.
\subsection{Estimating the update rate of the discount factor}

The update rule for the discount factor in (\ref{equ:rule2}) requires the access to the cost $J_{\gamma^k}(K^{k+1})$. To this end, we estimate $J_\gamma(K)$ by Monte Carlo sampling from the simulator as 
\begin{equation}\label{equ:Monte}
	\widehat{J}_{\gamma}^{\tau,N}(K) = \frac{1}{N}\sum_{i=1}^{N} V_{\gamma}^{\tau}(K, x_0^i),
\end{equation}
where $N\in \mathbb{N}_+$ and 
\begin{equation}\label{equ:costesti}
	\begin{aligned}
	&V_{\gamma}^{\tau}(K,x_0^i) = \sum_{t=0}^{\tau-1} (x_{t}^{\top} Q x_{t}+u_{t}^{\top} R u_{t})\\
	&\text {subject to} ~(\ref{equ:rescale_dyna})~ \text{and}~u_t = -K x_t,
	\end{aligned}
\end{equation}
where $x_0^i, ~i \in \{1,2,\dots,N\}$ are independently sampled initial states  from the distribution $\mathcal{D}$ and $\tau$ is the length of time horizon. We show that the estimate error $|\widehat{J}_{\gamma}^{\tau,N}(K)-J_{\gamma}(K)|$ induced by Monte Carlo sampling and truncation can be well controlled with high probability by using only a constant number of system rollouts $N$, the proof of which is provided in Appendix \ref{app:lemm5}.

\begin{lemma}\label{lem:error}
	For $K \in \mathcal{S}_{\gamma}$ and a constant $0<\delta<1$, let the horizon $\tau$ and the number of rollouts $N$ satisfy
	\begin{equation}\label{equ:condi}
		\tau \geq \frac{2J_{\gamma}(K)}{\underline{\sigma}(Q)}\log\left(\frac{J_{\gamma}(K)d^2}{\underline{\sigma}(Q)}\right), ~ N \geq 8d^4 \log\left(\frac{2}{\delta}\right).
	\end{equation}
	Then, with probability not smaller than $1-\delta$, it holds that
	\begin{equation}\label{equ:eval_bound}
	|J_{\gamma}(K) - \widehat{J}_{\gamma}^{\tau,N}(K) |\leq \frac{1}{2}J_{\gamma}(K).
	\end{equation}
\end{lemma}

To prove it, we show that an upper bound of $|V_{\gamma}^{\infty}(K, x_0) - V_{\gamma}^{\tau}(K, x_0)|$ decays exponentially with the horizon $\tau$ in Lemma \ref{lem:bounds_x}, which is a discrete-time counterpart of \cite[Lemma 11]{mohammadi2021convergence}. Then, we select the specific horizon and the accuracy to make the sample complexity $N$ independent of the discount factor. Finally, we use concentration bounds for the bounded random variable $x_0$ to derive (\ref{equ:condi}).

By Lemma \ref{lem:error}, the cost $J_{\gamma^k}(K^{k+1})$ in (\ref{equ:rule2}) can be upper bounded by $J_{\gamma^k}(K^{k+1}) \leq 2\widehat{J}_{\gamma^k}^{\tau,N}(K^{k+1})$ with high probability under the condition (\ref{equ:condi}). Thus, we can revise the update rate in (\ref{equ:rule2}) as  
\begin{equation}\label{equ:rule3}
\alpha^k = \frac{\underline{\sigma}(Q+(K^{k+1})^{\top}RK^{k+1})}{2\widehat{J}_{\gamma^k}^{\tau,N}(K^{k+1}) - \underline{\sigma}(Q+(K^{k+1})^{\top}RK^{k+1})}.
\end{equation}
Moreover, it follows from (\ref{equ:condi}) that the sample complexity of $N$ for estimating $\alpha^k$ remains constant in the iteration. In comparison of \cite{perdomo2021stabilizing}, it is impossible to derive such an explicit bound via the binary or random search.

Next, we adopt PG methods for the cost reduction step in Algorithm \ref{alg:exact}  via the simulator.

\subsection{Cost reduction via zero-order PG methods}
To achieve $J_{\gamma^k}(K^{k+1}) < D$ in (\ref{equ:cost_req}),  we use samples of the cost to estimate the gradient $\widehat{{\nabla}J_{\gamma}}(K)$ via zero-order methods, such as mini-batching~\cite{fazel2018global} and one- or two-point methods \cite{malik2019derivative, mohammadi2021convergence}. In Algorithm \ref{alg:gradient}, we adopt the two-point gradient estimate as it usually leads to a lower sample complexity for PG methods~\cite{malik2019derivative}. Then, the PG method with gradient estimate to solve \eqref{prob:dlqr} under a fixed $\gamma$ iterates as
\begin{equation}\label{equ:appro_grad}
K_{j+1} = K_j - \eta \widehat{{\nabla}J_{\gamma}}(K_j),
\end{equation}
where $\eta$ is the stepsize. In this section, we use the subscript $j$ to denote the iteration of the policy under a fixed discount factor and use the superscript $k$ to denote the iteration of the discount factor.
\begin{algorithm}[t]
	\caption{The two-point gradient estimate}
	\label{alg:gradient}
	\begin{algorithmic}[1]
		\Require
		A discount factor $\gamma$, a policy $K \in \mathcal{S}_{\gamma}$, a smoothing radius $r$, a time horizon $\tau$, a number of system rollouts $N_e$.
		\For{$i=1,\cdots,N_e$}
		\State Sample a perturbation matrix $U_{i}$ uniformly from the unit sphere $\mathcal{U}^{mn-1}$.
		\State Set $K_{i,1} = K + r\sqrt{mn}U_{i}$ and $K_{i,2} = K - r\sqrt{mn}U_{i}$.
		\State Sample an initial state $x_0^i$ using the distribution $\mathcal{D}$.
		\State Simulate system (\ref{equ:rescale_dyna}) to obtain $V_{\gamma}^{\tau}(K_{i,1}, x_0^i)$ and $V_{\gamma}^{\tau}(K_{i,2}, x_0^i)$.
		\EndFor	
		\Ensure
		Gradient estimate $$\widehat{\nabla J_{\gamma}}(K) = \frac{1}{2rN_e} \sum_{i=1}^{N_e}(V_{\gamma}^{\tau}(K_{i,1}, x_0^i) - V_{\gamma}^{\tau}(K_{i,2}, x_0^i))U_i.$$
	\end{algorithmic}
\end{algorithm}

\begin{algorithm}[t]
	\caption{The PG based explicit discount method}
	\label{alg:mfeda}
	\begin{algorithmic}[1]
		\Require
		An initial policy $K^0 = 0$, a discount factor $\gamma^0<1/\rho^2(A)$, an coefficient $0<\xi<1$, a number of PG iterations $M$.
		\For{$k=0,1,\cdots$}
		\State \textbf{Cost reduction:} Let $K_0=K^k$. Iterate
				\begin{equation}\label{equ:PGiter}
					K_{j+1} = K_j - \eta \widehat{{\nabla}J_{\gamma^k}}(K_j), ~j=0,1,\dots,M-1,
				\end{equation}	
		     ~~~	and set $K^{k+1} = K_{M}$.
		\State \textbf{Update of the discount factor:}  
		
		$\gamma^{k+1} = (1+\xi\alpha^k)\gamma^k$ with $\alpha^k$ given by (\ref{equ:rule3}).
		\If{$\gamma^{k+1} \geq  1 $} return a stabilizing policy $K^{k+1}$.
		\EndIf
		\EndFor	
	\end{algorithmic}
\end{algorithm}
Built upon Algorithm \ref{alg:exact}, we present the explicit discount method based on PG methods in Algorithm \ref{alg:mfeda}. We characterize how many PG iterations it needs to achieve (\ref{equ:cost_req}) and the corresponding sample complexity at the $k$-th iteration.

We first note that the PG update (\ref{equ:appro_grad}) meets global convergence under proper  stepsize, which is established by exploiting the gradient dominance property and local smoothness of the discounted LQR cost~\cite{mohammadi2021convergence}, i.e.,
\begin{align*}
&J_{\gamma}(K) - J_{\gamma}^* \leq \frac{1}{2\mu} \|\nabla J_{\gamma}(K)\|^2,~\text{and}~ \\
&J_{\gamma}(K') - J_{\gamma}(K) \leq \left<\nabla J_{\gamma}(K), K'-K\right> + \frac{L(a)}{2}\|K'-K\|_F^2,
\end{align*}
for all $K$ and $K'$ such that the line segment connecting them belongs to the sublevel set $\{K | J_{\gamma}(K) \leq a\}$, where $a$ is a positive constant, $\mu = {\underline{\sigma}(R)}/{(2\|\Sigma^*(\gamma)\|)}$ is the gradient dominance constant with $\Sigma^*(\gamma)$ being the solution of (\ref{def:cov}) under the optimal policy, and $L(a)$ is the smoothness parameter polynomial in the LQR parameters $(a, \|A\|,\|B\|,\|Q\|,\|R\|,\gamma)$. This means that the cost requirement (\ref{equ:cost_req}) can be ensured after a finite number of PG iterations, as shown below. 

\begin{lemma}\label{lem:grad}
	Consider the zero-order PG method in (\ref{equ:PGiter}) at the $k$-th iteration starting from $K_0 = K^k \in \mathcal{S}_{\gamma^k}$. Let $\widehat{{\nabla}J_{\gamma^k}}(K)$ be computed by Algorithm \ref{alg:gradient} where
	$\tau \geq \theta_1^k \log({1}/{(D - J_{\gamma^k}^*)}), N_e \geq c(1+\beta^{4} \phi^{4} \theta_2^k \log ^{6} n) n $, $r<\theta_3^k (D - J_{\gamma^k}^*)^{1/2}$, and $\eta \leq 1/(\omega^k L^k)$. 
	If 	
	\begin{equation}\label{equ:pgroll}
		j \geq -\log((J_{\gamma^k}(K^k)-J_{\gamma^k}^*)/(D - J_{\gamma^k}^*))/\log(1-\mu^k\eta/8),
	\end{equation}
	then $J_{\gamma^k}(K_j) \leq D$ with probability not smaller than
	\begin{equation}\label{equ:prob_pg}
		1-c' j(n^{-\beta}+N_e^{-\beta}+N_ee^{-\frac{n}{8}}+e^{-c' N_e}).
	\end{equation}
	 Here, $\omega^k = c''(\sqrt{m}+\beta \phi^{2} \theta_2^k \sqrt{m n} \log n)^{2}$, $\phi$ is a function of $d$, $\mu^k$ is the gradient dominance constant, $L^k$ is the smoothness parameter over the sublevel set $\{K | J_{\gamma^k}(K) \leq J_{\gamma^k}(K^k)\}$, $\beta, c,c',c''$ are positive constants, and $L^k, \theta_1^k, \theta_2^k, \theta_3^k$ are polynomials in the parameters of the discounted LQR problem $(J_{\gamma^k}(K^k), \|A\|,\|B\|,\|Q\|,\|R\|, \underline{\sigma}(Q), \underline{\sigma}(R), \gamma^k)$.
\end{lemma}

The proof is provided in Appendix \ref{app:grad} and follows from  \cite[Theorem 1]{mohammadi2020linear} by focusing on our discounted problem and the cost requirement (\ref{equ:cost_req}).

\subsection{Convergence and sample complexity}

We now show the convergence and sample complexity result of Algorithm \ref{alg:mfeda}. 
{As many variants of LQR problems in \cite{zhang2021policy,zhao2021global, zhang2019policy, zheng2021analysis, duan2022optimization,furieri2020learning, li2021distributed, fatkhullin2021optimizing,  gravell2020learn, tu2019gap, yang2019provably}, we also adopted the convergence analysis approaches in   \cite{fazel2018global, bu2019lqr, malik2019derivative, mohammadi2020linear, mohammadi2021convergence} in Section IV. Since each LQ problem has its unique features, using those approaches to a specific LQ problem is non-trivial.}
{Specifically, here we are dealing with a sequence of discounted LQR problems where both the discount factor and the initial cost for each one are changing, and require to provide uniform bounds for the increasing LQR costs. Moreover, algorithm parameters need to be carefully designed to study the convergence and sample complexity. Thus, the existing references cannot be used to solve these technical challenges. }

{The sample complexity of Algorithm \ref{alg:mfeda} is formally given below, the proof of which is provided in Appendix \ref{app:main}.}
{\begin{theorem}\label{thm:sample}
Let the parameters in  Algorithm \ref{alg:mfeda} satisfy
\begin{equation}\label{condition}
	\begin{aligned}
	& \tau \geq \max\left\{\frac{4J^*_1}{\underline{\sigma}(Q)}\log\left(\frac{2J_1^*d^2}{\underline{\sigma}(Q)}\right),  \overline{\theta}_1 \log\left(\frac{1}{J_1^*}\right)\right\}, \\
	& N \geq 8d^4 \log({2}/{\delta}), N_e \geq c(1+\beta^{4} \phi^{4}  \overline{\theta}_2 \log ^{6} n) n, \\
	& r<\underline{\theta}_3 \sqrt{J_1^*}, \eta \leq \frac{1}{\overline{w}\overline{L}}, 
	 M \geq \frac{\overline{L} J^*_1\overline{\omega}}{\underline{\sigma}(R)} \log \left(\frac{12J^*_1}{(1-\xi)\underline{\sigma}(Q)}\right).
	\end{aligned}
\end{equation}	
	If $k \geq {(6J_1^*-\underline{\sigma}(Q))}\log({1}/{\gamma^0})/(\xi\underline{\sigma}(Q))$, then with probability not smaller than 
	\begin{equation}\label{equ:prob}
		1 - k\left(\delta + c' M(n^{-\beta}+N_e^{-\beta}+N_ee^{-\frac{n}{8}}+e^{-c' N_e})\right),
	\end{equation}
	Algorithm \ref{alg:mfeda} returns a stabilizing policy $K^k$ with \begin{equation}\label{equ:perf}
		\rho(A-BK^k)<\sqrt{1-(1-\xi)\underline{\sigma}(Q)/(6J_1^*)}.
	\end{equation} 
	Here, $\overline{w} = c''(\sqrt{m}+\beta \phi^{2} \overline{\theta}_2 \sqrt{m n} \log n)^{2}$, $\overline{\theta}_1, \overline{\theta}_2, \underline{\theta}_3, \overline{L}$ are positive polynomials in $(\|A\|,\|B\|,\|Q\|,\|R\|,\underline{\sigma}(Q), \underline{\sigma}(R))$, $0< \delta <1$ is a user-specified constant, and $\phi,\beta, c,c',c''$ are the constants in Lemma \ref{lem:grad}.
\end{theorem}}

To prove it, we first derive uniform bounds (\ref{condition}) for the algorithm parameters in Lemma \ref{lem:error} and Lemma \ref{lem:grad}, and then use union bounds to provide the probability of convergence in~\eqref{equ:prob}. Particularly, we are able to provide a lower bound for $M$ independent of $k$ thanks to our explicit update rule. In contrast, Ref. \cite{perdomo2021stabilizing} cannot provide such a bound due to the binary or random search of the discount factor.




The total number of system rollouts required  to find a stabilizing policy under the probability (\ref{equ:prob}) is 
\begin{align*}
& k(N + N_eM) \approx \\
&\frac{6J^*_1-\underline{\sigma}(Q)}{\xi\underline{\sigma}(Q)}\log\left(\frac{1}{\gamma^0}\right) \bigg(8d^4 \log\left(\frac{2}{\delta}\right) +  (1+\beta^{4} \overline{\theta}_2^2\phi^{4}  \log^{6} n)  \\
&\times \frac{CJ^*_1}{\underline{\sigma}(R)} (1+\beta \overline{\theta}_2 \phi^{2} \sqrt{n} \log n)^{2} mn\log \left(\frac{12J^*_1}{(1-\xi)\underline{\sigma}(Q)}\right) \bigg),
\end{align*}
where $C$ is a constant. To show the dependence on important system parameters $n,m,\rho(A)$, we let $\gamma^0 = 1/(2\rho^2(A))$ and write the sample complexity more compactly as 
\begin{equation}\label{equ:sam}
	\log(\rho(A))\cdot \widetilde{\mathcal{O}}(n^2m)\cdot\text{poly}(\|A\|,\|B\|,\|Q\|,\|R\|,J^*).
\end{equation}
This result reveals that the sample complexity is quadratic in the system dimension $n$ (by neglecting the coefficient polynomial in $\log(n)$) and linear in the input dimension $m$, which matches the state-of-the-art complexity of PG methods for solving the LQR problem \cite{mohammadi2020linear,mohammadi2021convergence}. In fact, it amounts to multiplying a logarithmic coefficient $\log(\rho(A))$ to that for solving the LQR problem \cite{mohammadi2020linear,mohammadi2021convergence}. This means that the difficulty of stabilizing an unstable linear system increases with the maximal unstable eigenvalue of the system. The complexity also has a slight dependence on the user-specified coefficient $\xi$, which is designed to show the performance guarantees of the resulted policy in \eqref{equ:perf}. 

Our sample complexity quadratic in $n$ is comparable with that of the model-based approach, which has a dimension dependence linear in $n$~\cite{malik2019derivative}. This is because each observation in model-based approach is a state vector with dimension $n$, while herein a sample of the cost is a scalar. 

Finally, we remark some extensions of our explicit discount method. First, we can leverage the explicit update rule to provide finite-time convergence guarantees and sample complexity for other RL based discount methods,  e.g.,~\cite{jing2021learning,lamperski2020computing}. {Second, our PG based discount method is straightforward to extend to linear systems with stochastic noises by leveraging the equivalence between LQR formulations with noises and random initial conditions established in \cite[Lemma 10]{malik2019derivative}. Third, it can be used to stabilize nonlinear systems around the equilibrium point, which we shall confirm in the simulation.} 



\section{Simulation}\label{sec:exp}

In this section, we first validate the theoretical results of our PG based discount method in linear systems. Then, we perform it over a classical nonlinear system to show the effectiveness of the proposed method. We observe that our method results in a larger Region of Attraction (ROA) of initial states than that of using the LQR of the system linearized at the equilibrium point.

\subsection{Linear systems}

We first perform simulations over a two-dimensional example for illustration. Consider the following unstable dynamical system with single control input 
$$
A = \begin{bmatrix}
4 &3 \\
3 &1.5
\end{bmatrix}, ~~ B = \begin{bmatrix}
2 \\
2
\end{bmatrix}, ~~Q = \begin{bmatrix}
1 & 0\\
0 & 1
\end{bmatrix}, ~~R = 2.
$$
Clearly, $(A,B)$ is controllable. Let the initial state distribution $\mathcal{D}$ be the standard normal distribution. The selections of the parameters in Algorithm \ref{alg:mfeda} are as follows. The number of rollouts in function evaluation is set to $N = 20$, and the horizon is set to $\tau = 100$. We select an initial policy $K^0 = 0$, the coefficient $\xi = 0.9$ and discount factor $\gamma^0 = 10^{-3} < 1/\rho^2(A) = 1/36$. We apply one-step gradient descent each iteration in Algorithm \ref{alg:mfeda}, i.e., the policy is updated only once for each $k$. The stepsize of gradient descent is set to a constant $\eta = 10^{-3}$, and the smooth radius and the number of rollouts for gradient estimation in Algorithm \ref{alg:gradient} are set to $r = 2\times 10^{-3}$ and $N_e =20$, respectively. Fig. \ref{pic:pgland} illustrates the policy trajectory. In less than $100$ iterations, Algorithm \ref{alg:mfeda} returns a stabilizing policy. 
\begin{figure}[t]
	\centering
	\subfigure[Convergence of Algorithm \ref{alg:mfeda}.]{
		\includegraphics[height=30mm]{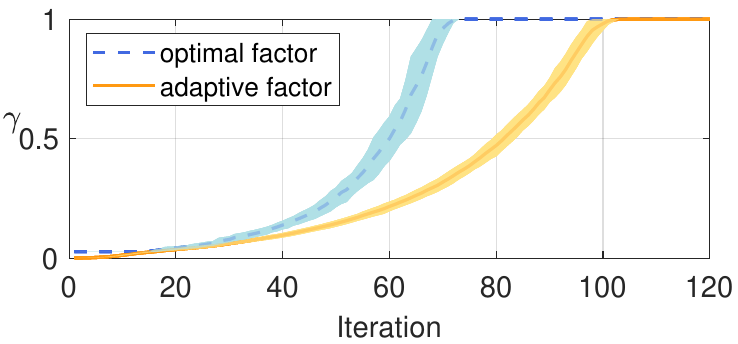}}
	\subfigure[Model-based implementation.]{
		\includegraphics[height=30mm]{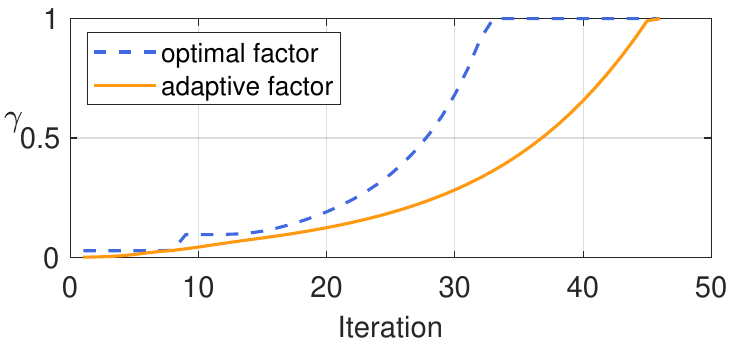}}
	\caption{Convergence of the discount factor. The centreline denotes the mean of 20 independent trials and the shaded region demonstrates their standard deviation.}
	\label{pic:discount}
\end{figure}

Next, we show the convergence of the discount factor in Fig \ref{pic:discount}, where we report the results of $20$ independent trails. From $\gamma^0 = 10^{-3}$, the adaptive discount factor (yellow solid line) grows almost exponentially to $1$ within $100$ iterations. The sample complexity can be calculated by $100\times(N+N_e)$, i.e.,  total number of $4\times 10^3$ rollouts with horizon $\tau = 100$. For comparison, we also plot the maximal discount factor under $K^k$ each iteration, given as $\gamma^k_{\text{opt}} = 1/\rho^2(A-BK^k)$, denoted by the ``optimal'' discount factor (blue dashed line). We observe that the update rate $\alpha^k$ in (\ref{equ:rule3}) well approximates the upper bound in the first $50$ iterations, and the gap increases reasonably due to the approximation of condition (\ref{equ:rule}). Moreover, the variance induced by independent trials is competitively small, considering our low sample complexity. We also display a model-based implementation of Algorithm \ref{alg:mfeda}, where we assume that $(A,B)$ is known and compute the update rate  by (\ref{equ:coro}) in Theorem \ref{coro}. It is shown that using the exact $J_{\gamma}(K)$ and the gradient, the required number of iterations can be reduced to less than $50$. 

\begin{figure}[t]
	\centering
	\includegraphics[height=45mm]{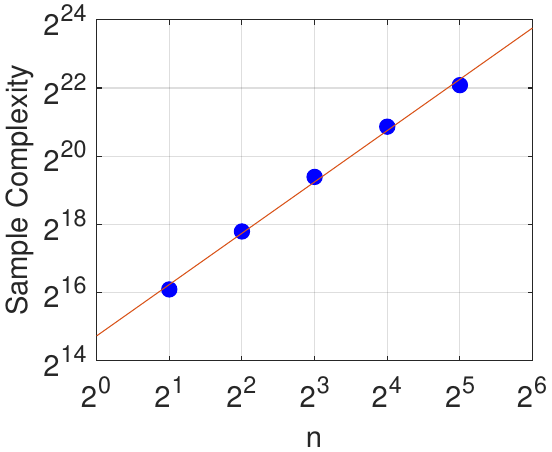}
	\caption{Total number of system rollouts as a function of state dimension $n$.
	}
	\label{pic:n}
\end{figure}

\begin{figure}[t]
	\centering
	\includegraphics[height=45mm]{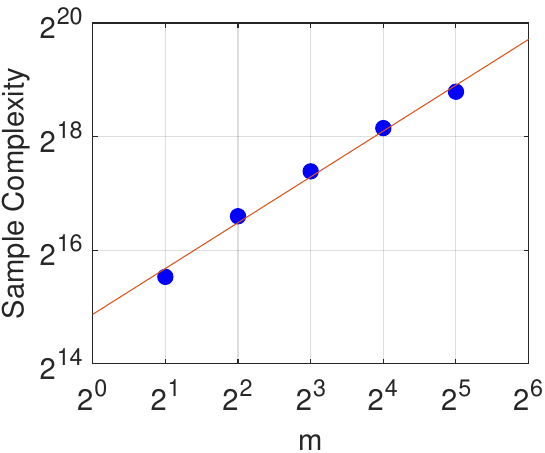}
	\caption{Total number of system rollouts as a function of input dimension $m$.
	}
	\label{pic:m}
\end{figure}

To validate the sample complexity result in Theorem \ref{thm:sample}, we study the total number of system rollouts required by Algorithm \ref{alg:mfeda} as a function of the state dimension $n$ and the input dimension $m$. We first validate the dependence on $n$ in (\ref{equ:sam}), where we fix the system parameters $m=8$, $Q$ and $R$ as identity matrices, and randomly sample $A = 2(\tilde{A}+\tilde{A}^{\top})/\|\tilde{A}+\tilde{A}^{\top}\|$ with each element of $\tilde{A}$ subject to the standard normal distribution and $B=\tilde{B}/\|\tilde{B}\|$ with $\tilde{B}$ also being Gaussian such that $\|A\|=2,\|B\|=1$ and $(A,B)$ is stabilizable with probability one. For algorithm parameters, we set the initial discount factor to $0.9/\rho^2(A) = 0.225$, $N_e = 20\times n$ as suggested by Lemma \ref{lem:grad}, and other parameters the same as before. For each $n$, we repeat $20$ independent trials and the mean of the sample complexity is recorded. The result is plotted in Fig. \ref{pic:n}. It can be observed that the sample complexity is approximately quadratic in the state dimension $n$. Then, we validate the dependence on $n$ in (\ref{equ:sam}), where $n=2$ is fixed and other parameters are set as before. The result is plotted in Fig. \ref{pic:m}. As indicated by (\ref{equ:sam}), the sample complexity is approximately linear in the input dimension $m$.

\subsection{Nonlinear systems}
In this subsection, we perform our explicit discount method over the classical cart-pole system with the following continuous dynamics \cite{perdomo2021stabilizing}
$$
\begin{bmatrix}
m_{p}+m_{c} & -m_{p} l \cos (\theta) \\
-m_{p} l \cos (\theta) & m_{p} l^{2}
\end{bmatrix}
\begin{bmatrix}
\ddot{z} \\
\ddot{\theta}
\end{bmatrix}=\begin{bmatrix}
u-m_{p} l \sin (\theta) \dot{\theta}^{2} \\
m_{p} g l \sin (\theta)
\end{bmatrix},
$$
where $z$ and $\dot{z}$ denote the position and velocity of the cart, $\theta$ and $\dot{\theta}$ denote the angular position and velocity of the pole, and the control input $u$ is the horizontal force to the cart. The state is defined as $x= (z,\theta, \dot{z},\dot{\theta})$. The goal is to regulate the system state to attain the upright equilibrium with $\theta = 0$. The system parameters include the mass of the pole $m_p$ and the cart $m_c$, the length of the pendulum $l$, and the gravity acceleration $g$. In the simulation, we set them to $m_p = m_c = l= g =1$.

For implementation, we use the forward Euler rule to discretize the continuous dynamics with sampling time $dt = 0.02s$. We limit the policy to be linear feedback $u = -Kx$ and implement Algorithm \ref{alg:mfeda} to search a stabilizing gain over a neighbourhood around the equilibrium point. We set the penalty matrices to $Q = 2\times I_{4\times4}$ and $R= 1$. The initial discount factor is $\gamma^0 = 0.01$. Other parameters are set to $\tau=10^3, N_e=N=20, r=0.01, \eta=10^{-3}, \xi=1$. In the simulator, we obtain samples of the following cost
\begin{align*}
&J_{\gamma}^{\text{nl}}(K) = \mathbb{E}_{x_0} \sum_{t=0}^{\infty} (x_{t}^{\top} Q x_{t}+u_{t}^{\top} R u_{t}),\\
&\text {subject to} ~x_{t+1} = \sqrt{\gamma}f(x_t,u_t), x_0\sim \mathcal{D}_{\text{nl}}  ~ \text{and}~~u_t = -K x_t,
\end{align*}
where $x_{t+1} = f(x_t,u_t)$ denotes the discrete-time cart-pole system and the distribution $\mathcal{D}_{\text{nl}}$ of the initial state $x_0$ is uniform over the cube $[-r_{\text{ini}}, r_{\text{ini}}]^4$.
Since the covariance of $\mathcal{D}_{\text{nl}}$ is not unitary, the update rate at the $k$-th iteration is accordingly modified as 
$$
\alpha^k = \frac{\underline{\sigma}(Q+(K^{k+1})^{\top}RK^{k+1})}{3J_{\gamma^k}^{\text{nl}}(K^{k+1})/r_{\text{ini}} - \underline{\sigma}(Q+(K^{k+1})^{\top}RK^{k+1})},
$$
where $J_{\gamma^k}^{\text{nl}}(K^{k+1})$ is estimated from the simulator.

To evaluate the performance of the resulted policy, we introduce the notion of ROA which describes the region of the initial states that can asymptotically converge to the equilibrium point. Overall, it reflects the ability of a controller to locally stabilize the system. In this simulation, we use the largest ball of the initial states that asymptotically converge to the origin to estimate the ROA based on Monte Carlo sampling. To this end, we randomly select an initial state from the sphere $r_{\text{test}}\cdot\mathcal{U}^{3}$ and simulate the closed-loop nonlinear dynamics for $2\times 10^3$ steps. We repeat it for $10^3$ times to see whether $\|x_{2\times 10^3}\| \leq r_{\text{test}}/5$ (the state converges) in all cases. If the answer is positive, then we conclude that $r_{\text{roa}} > r_{\text{test}}$. We then increase $r_{\text{test}}$ and repeat the test until for a radius $\bar{r}_{\text{test}}$ the state norm fails to converge to within $\bar{r}_{\text{test}}/5$. Hence, we can estimate the ROA as the ball with radius $r_{\text{roa}} = \bar{r}_{\text{test}}$.
\begin{table}[!t]\label{table}
	\caption{The region of attraction radius $r_{\text{roa}}$.}
	\begin{center}
		\begin{tabular}{|ccccc|}
			\hline
			$r_{\text{ini}}$ & $0.1$ & $0.3$ & $0.5$ & LQR\\
			\hline
			$r_{\text{roa}}$ & $0.61$ & $0.70$ & $0.79$& $0.66$\\
			\hline
		\end{tabular}
	\end{center}
\end{table}

We perform our algorithm with initial state radius $r_{\text{ini}} = 0.1,0.3,0.5$, each with $3$ independent trails, and estimate the average ROA radius of the resulting controllers. We also compare the results with the optimal LQR controller of the linearized dynamical model at the origin. Our result is summarized in Table I. Interestingly, though the LQR controller has the well-known stability margin, our explicit discount method algorithm returns a stabilizing policy with a larger ROA radius when $r_{\text{ini}} = 0.3,0.5$. We also observe that the ROA radius becomes slightly larger as $r_{\text{ini}}$ increases. This is due to that a large $r_{\text{ini}}$ can facilitate the exploration of the unknown nonlinear system, which helps learn a better controller.

\section{Conclusion and future work}\label{sec:conclusion}

In this paper, we have proposed a PG method to learn a stabilizing policy for linear systems by designing an update rule for the discount factor. Moreover, we have shown that the algorithm converges in finite iterations and proved explicit sample complexity in terms of the total number of system rollouts. Simulations have validated the theoretical results and shown the effectiveness of the proposed method on nonlinear systems. As the first to provide a formal sample complexity guarantee for the stabilization problem, we hope our work helps pave the way for understanding the performance of PG methods for general systems. 

We now discuss some future directions. The prior work \cite{perdomo2021stabilizing} has designed a binary search procedure of the noisy cost to update the discount factor, which is applicable to smooth nonlinear dynamics. Since our update rule for the discount factor is derived based on the structure of linear systems, it would be valuable to examine whether it can also be extended to nonlinear systems from a theoretical perspective. It would also be interesting to further exploit the connections between the notion in \cite{safonov1977gain} and ours of the stability margin, as noted by Remark \ref{remark1}. The cost reduction step in Algorithm \ref{alg:exact} can be replaced with other RL methods (e.g., Q-learning and actor-critic) to study their performance and sample complexities for stabilization. Another direction is to study the PG methods for finite-horizon LQR problems with a single feedback gain $K$. Similar to the discount methods, the idea is to gradually increase the horizon to find a stabilizing gain. The major challenge is that the optimal policy for finite-horizon LQR may not be unique and the optimization landscape is unclear yet. Despite solving the stabilization problem, we believe it is valuable to further investigate the role of the discount factor in the PG methods of LQR problems. Existing literature in the RL field has shown that the discount factor works as a regularizer for generalization~\cite{amit2020discount}, and can even accelerate the convergence~\cite{franccois2015discount}. However, a systematic understanding of such phenomena is still lacking in the LQR setting, which is left as our future work.

\section*{Acknowledgement}
We would like to sincerely thank the associate editor and anonymous reviewers for their constructive comments, which significantly improve the quality of this article.

\bibliographystyle{IEEEtran}
\bibliography{mybibfile}

\appendices
\section{Proof of Theorem \ref{thm:tech}}\label{app:Lemma 5}
We first show that the update rate has a uniform lower bound, i.e., 
\begin{equation}\label{equ:dfs}
	\alpha^k \geq {\underline{\sigma}(Q)}/{(D-\underline{\sigma}(Q))}, \forall k\in \mathbb{N}.
\end{equation}
By (\ref{equ:cost_req}), the cost reduction step returns a policy such that $J_{\gamma^{k}}(K^{k+1}) \leq D$. Thus, $\alpha^k$ is lower bounded by 
$$
\alpha^k = \frac{\underline{\sigma}(Q+(K^{k+1})^{\top}RK^{k+1})}{J_{\gamma^k}(K^{k+1}) - \underline{\sigma}(Q+(K^{k+1})^{\top}RK^{k+1})} \geq  \frac{\underline{\sigma}(Q)}{D - \underline{\sigma}(Q)}.
$$

Then, we show that $J_{\gamma^k}(K^{k})$ has a uniform upper bound \begin{equation}\label{equ:asb}
	J_{\gamma^k}(K^{k}) \leq {D^2}/{((1-\xi)\underline{\sigma}(Q))}, \forall k\in \mathbb{N}_+.
\end{equation}
To this end, we show that the spectral radius of the damped system is upper bounded, 
\begin{equation}\label{equ:spectralbound}
\sqrt{\gamma^{k}}\rho(A-BK^{k}) \leq \sqrt{1-{(1-\xi)\underline{\sigma}(Q)}/{D}}, \forall k\in \mathbb{N}_+.
\end{equation}

Since at the $k$-th iteration $\sqrt{\gamma^k}\rho(A-BK^{k+1})<1$, Theorem \ref{coro} implies that $\sqrt{(1+\alpha^k)\gamma^{k}}(A-BK^{k+1})$ is stable with $\alpha^k$ given by (\ref{equ:rule2}). Then, the new discount factor $\gamma^{k+1} = (1+\xi\alpha^k)\gamma^k$ yields that
\begin{align*}
&\sqrt{\gamma^{k+1}}\rho(A-BK^{k+1})=\sqrt{(1+\xi\alpha^k)\gamma^k}\rho(A-BK^{k+1})\\
&= \frac{\sqrt{(1+\xi\alpha^k)\gamma^k}}{\sqrt{(1+\alpha^k)\gamma^k}}\cdot\sqrt{(1+\alpha^k)\gamma^k}\rho(A-BK^{k+1})\\
&< \frac{\sqrt{(1+\xi\alpha^k)\gamma^k}}{\sqrt{(1+\alpha^k)\gamma^k}} = \sqrt{1 - \frac{1-\xi}{(\alpha^k)^{-1}+1}}.
\end{align*}

By (\ref{equ:dfs}), it further leads to 
$$
\sqrt{\gamma^{k+1}}\rho(A-BK^{k+1}) \leq \sqrt{1-{(1-\xi)\underline{\sigma}(Q)}/{D}}.
$$	

Then, we prove that $J_{\gamma^{k+1}}(K^{k+1})$ has a uniform upper bound. Since the PG descent step in line 2 ensures $J_{\gamma^k}(K^{k+1})\leq D$, it follows from the trace inequality that
\begin{align*}
D &\geq J_{\gamma^k}(K^{k+1})\\
& = \text{Tr}((Q+(K^{k+1})^{\top}RK^{k+1})\Sigma_{K^{k+1}}(\gamma^k)) \\
&\geq \text{Tr}((Q+(K^{k+1})^{\top}RK^{k+1}) \underline{\sigma}( \Sigma_{K^{k+1}}(\gamma^k))) \\
&\geq \text{Tr}(Q+(K^{k+1})^{\top}RK^{k+1}),
\end{align*}
where the last inequality follows from Lemma \ref{lem:closed} that $\Sigma_{K^{k+1}}(\gamma^k) \geq I$. Jointly with  Lemma \ref{lem:closed}, it implies that
\begin{equation}\label{equ:jj}
\begin{aligned}
J_{\gamma^{k+1}}(K^{k+1}) &= \text{Tr}((Q+(K^{k+1})^{\top}RK^{k+1})\Sigma_{K^{k+1}}(\gamma^{k+1}))\\
&\leq \text{Tr}(Q+(K^{k+1})^{\top}RK^{k+1})\|\Sigma_{K^{k+1}}(\gamma^{k+1})\| \\
&\leq D\|\Sigma_{K^{k+1}}(\gamma^{k+1})\|.
\end{aligned}
\end{equation}

Let $\Sigma^k := \Sigma_{K^{k}}(\gamma^{k})$ for short. By iterating (\ref{def:cov}), it holds that $\Sigma^k = \sum_{t=0}^{\infty}(\gamma^k)^t((A-BK^k)^t)^{\top}(A-BK^k)^t$. It follows from (\ref{equ:spectralbound})  that 
$$
\|\Sigma^k\| \leq  \sum_{t=0}^{\infty}(1-\frac{(1-\xi)\underline{\sigma}(Q)}{D})^{t} = \frac{D}{(1-\xi)\underline{\sigma}(Q)}, \forall k\in \mathbb{N}_+.
$$
Inserting the above into (\ref{equ:jj}) yields (\ref{equ:asb}).

Finally, it follows from (\ref{equ:dfs}) that Algorithm \ref{alg:exact} halts within $${\log(1/\gamma^0)}/({\log(1+{\xi\underline{\sigma}(Q)}/{(D - \underline{\sigma}(Q))})})$$ iterations. Under the approximation $\log(1+x)\approx x$, the above is simplified as  ${(D-\underline{\sigma}(Q))}\log ({1}/{\gamma^0})/({\xi\underline{\sigma}(Q)})$. Jointly with  (\ref{equ:asb}), Algorithm \ref{alg:exact} returns a stabilizing policy with the spectral radius of the closed-loop system being upper bounded by $\sqrt{1-{(1-\xi)\underline{\sigma}(Q)}/{D}}$.

\section{Proof of Lemma \ref{lem:error}}\label{app:lemm5}

We first present a technical lemma, which is a discrete-time counterpart of \cite[Lemma 11]{mohammadi2021convergence}.
\begin{lemma}\label{lem:bounds_x}
 Consider the system (\ref{equ:rescale_dyna}).	Let $K\in \mathcal{S}_{\gamma}$ and $P$ satisfy (\ref{equ:Pi}). Then, for $t \in \mathbb{N}$, it holds that
	$$
	\|x_t\|^2 \leq \left(1 - \frac{\underline{\sigma}(Q)}{\|P\|}\right)^t \frac{\text{Tr}(P)}{\underline{\sigma}(P)} \|x_0\|^2.
	$$
\end{lemma}
\begin{proof}
	The function $V(x) = x^{\top}Px$ is a Lyapunov function of the system $x_{t+1} = \sqrt{\gamma}(A-BK)x_t$ since 
	$$
	V(x_{t+1}) - V(x_t) \leq -c V(x_t)
	$$
	with $c:={\underline{\sigma}(Q)}/{\|P\|}$. Then, from an initial state $x_0$ we have
	$$
	 \underline{\sigma}(P)\|x_t\|^2 \leq V(x_t) \leq (1 - c )^tV(x_0)\leq (1 - c )^t \|x_0\|^2\text{Tr}(P),
	$$
	which completes the proof.
\end{proof}

We first establish an upper bound of the bias $|V_{\gamma}^{\infty}(K, x_0) - V_{\gamma}^{\tau}(K, x_0)|$ induced by finite horizon, which has exponential dependence on $\tau$. For a policy  $K\in \mathcal{S}_{\gamma}$, it follows that
\begin{align*}
&V_{\gamma}^{\infty}(K, x_0) - V_{\gamma}^{\tau}(K, x_0) \\
&=  \sum_{t=0}^{\infty} (x_{t}^{\top} Q x_{t}+u_{t}^{\top} R u_{t}) -\sum_{t=0}^{\tau -1} (x_{t}^{\top} Q x_{t}+u_{t}^{\top} R u_{t})  \\
&= \sum_{t=\tau}^{\infty} (x_{t}^{\top} Q x_{t}+u_{t}^{\top} R u_{t})= x_{\tau}^{\top}Px_{\tau} \leq \text{Tr}(P)\|x_\tau\|^2,
\end{align*}
where $P$ is the solution of (\ref{equ:Pi}) and $x_{\tau} = (\sqrt{\gamma}(A-BK))^{\tau}x_0$. 

By Lemma \ref{lem:bounds_x}, the above is further upper bounded by
\begin{align*}
	&V_{\gamma}^{\infty}(K, x_0) - V_{\gamma}^{\tau}(K, x_0)
	\leq  \left(1 - \frac{\underline{\sigma}(Q)}{\|P\|}\right)^t \frac{(\text{Tr}(P))^2}{\underline{\sigma}(P)} \|x_0\|^2 \\
	&\leq  \frac{d^2J_{\gamma}^2(K)}{\underline{\sigma}(Q)} \left(1 - \frac{\underline{\sigma}(Q)}{J_{\gamma}(K)}\right)^t.
\end{align*}

Then, we use concentration inequalities to bound the total error $|J_{\gamma}(K) - \widehat{J}_{\gamma}^{\tau,N}(K)|$. Let $x_0^i, i\in \{1,2,\dots,N\}$ be $N$ random initial states sampled independently from $\mathcal{D}$. 
Since the support of $\mathcal{D}$ is bounded by $\|x_0\| \leq d$, the random variable $V_{\gamma}^{\infty}(K, x_0^i)$ is bounded by $ 0\leq V_{\gamma}^{\infty}(K, x_0^i) = \text{Tr}(P x_0x_0^{\top}) \leq J_{\gamma}(K)d^2, \forall i \in \{1,2,\dots,N\}$. For a given constant $\epsilon$, we let the horizon be 
$$\tau = \frac{-\log(2J_{\gamma}^2(K)d^2/\underline{\sigma}(Q)\epsilon)}{\log(1- \underline{\sigma}(Q)/J_{\gamma}(K) )} \approx \frac{J_{\gamma}(K)}{\underline{\sigma}(Q)}\log\left(\frac{2J_{\gamma}^2(K)d^2}{\underline{\sigma}(Q)\epsilon}\right)
$$
such that $V_{\gamma}^{\infty}(K, x_0^i) - V_{\gamma}^{\tau}(K, x_0^i) \leq \epsilon/2$. Then, the following holds with probability $1$
\begin{equation}\label{equ:error}
\left|\frac{1}{N}\sum_{i=0}^{N-1} V_{\gamma}^{\tau}(K, x_0^i) -\frac{1}{N}\sum_{i=0}^{N-1} V_{\gamma}^{\infty}(K, x_0^i)\right| \leq \frac{\epsilon}{2}.
\end{equation}
Hence, the Hoeffding's inequality for the random variable $V_{\gamma}^{\infty}(K,x_0)$ yields that
\begin{align*}
&\text{Pr}(|\widehat{J}_{\gamma}^{\tau,N}(K) - J_{\gamma}(K)| \leq \epsilon)\\
&= \text{Pr}\left(|\frac{1}{N}\sum_{i=0}^{N-1} V_{\gamma}^{\tau}(K, x_0^i) - \mathbb{E}_{x_0}[V_{\gamma}^{\infty}(K,x_0^i)]| \leq \epsilon\right) \\
&\geq \text{Pr}\left(|\frac{1}{N}\sum_{i=0}^{N-1} V_{\gamma}^{\infty}(K, x_0^i) - \mathbb{E}_{x_0}[V_{\gamma}^{\infty}(K,x_0^i)]| \leq \frac{\epsilon}{2}\right) \\
& \geq 1- 2\exp\left(-\frac{N\epsilon^2}{2J_{\gamma}(K)^2d^4}\right),
\end{align*}
where the first inequality follows from (\ref{equ:error}) and the triangle inequality. Let $\delta = 2\exp(-{N\epsilon^2}/(2J_{\gamma}(K)^2d^4))$. Then, we conclude that at least with probability $1-\delta$, the estimation error is bounded by $|\widehat{J}_{\gamma}^{\tau,N}(K) - J_{\gamma}(K)| \leq \epsilon$. Letting $\epsilon = J_{\gamma}(K)/2$ completes the proof.

\section{Proof of Lemma \ref{lem:grad}}\label{app:grad}
Due to the equivalence between the discounted LQR in (\ref{prob:dlqr}) and the undiscounted LQR with damped systems in (\ref{equ:equivalent}), we can leverage the results \cite[Theorem 1]{mohammadi2020linear} for standard LQR to our problem. Note that the distribution of their initial state is assumed to have a bounded sub-Gaussian norm. Since the distribution $\mathcal{D}$ has bounded supports, there exists $\phi>0$ as a function of $d$ such that $\|x_0\|_{\Psi_2} \leq \phi$. Since we require that $J_{\gamma^k}(K_j) < D$, the optimality gap is no more than $\epsilon=D-J^*_{\gamma^k}$. Then, noting that our PG update starts from $K^k$ over the sublevel set $\{K | J_{\gamma}(K) \leq J_{\gamma}(K^k)\}$ and applying \cite[Theorem 1]{mohammadi2020linear} complete the proof.

\section{Proof of Theorem \ref{thm:sample}}\label{app:main}
The proof of Theorem \ref{thm:sample} relies on the following lemma.

\begin{lemma}\label{lem:newtech}
	Under the condition \eqref{condition}, the update rate has a uniform lower bound 
	\begin{equation}\label{equ:updae}
	\alpha^k \geq {\underline{\sigma}(Q)}/{(6J_1^*-\underline{\sigma}(Q))}
	\end{equation}
	with probability not smaller than 
	$$
	1 - k\left(\delta + c' M(n^{-\beta}+N_e^{-\beta}+N_ee^{-\frac{n}{8}}+e^{-c' N_e})\right), \forall k\in \mathbb{N}.
	$$ 
\end{lemma}
\begin{proof}
	Let $D = 2J_1^*$. We first show at the $k$-th iteration that if the PG updates successfully return $J^k(K^{k+1}) \leq D$ and the cost estimate satisfies $|J_{\gamma^k}(K^{k+1}) -\widehat{J}_{\gamma^k}^{\tau,N}(K^{k+1})| \leq J_{\gamma^k}(K^{k+1})/2$, then \eqref{equ:updae} holds with probability $1$. Then, we derive their failure probabilities and use union bound to complete the proof.
	
	Since $J_{\gamma^k}(K^{k+1}) \leq 2\widehat{J}_{\gamma^k}^{\tau,N}(K^{k+1})$, the update rate $\alpha^k$ in (\ref{equ:rule3}) satisfies Theorem \ref{coro}, which is lower bounded by 
	$$
	\alpha^k \geq \frac{\underline{\sigma}(Q)}{2\widehat{J}_{\gamma^k}^{\tau,N}(K^{k+1}) - \underline{\sigma}(Q)} \geq \frac{\underline{\sigma}(Q)}{3D - \underline{\sigma}(Q)} = \frac{\underline{\sigma}(Q)}{6J_1^* - \underline{\sigma}(Q)}.
	$$
	
	Now, we derive the failure probability. It is straightforward to check that the parameters in (\ref{condition}) satisfy the assumptions in Lemma \ref{lem:error} and Lemma \ref{lem:grad}. Clearly, $\tau$ and $N$ in (\ref{condition}) are sufficient conditions for (\ref{equ:condi}) since $J_{\gamma^k}(K^{k+1})\leq D = 2J_1^*$. Thus, at the $k$-th iteration, $|J_{\gamma^k}(K^{k+1}) -\widehat{J}_{\gamma^k}^{\tau,N}(K^{k+1})| \leq J_{\gamma^k}(K^{k+1})/2$ has the failure probability not larger than $\delta$. 
	
	Similarly, we next show that $\tau,N_e,r,\eta, M$ in (\ref{condition}) satisfy the condition in Lemma \ref{lem:grad}. We first note that $J^{k}(K^{k}) \leq {12(J_1^*)^2}/{((1-\xi)\underline{\sigma}(Q))}$, which follows from the derivation of (\ref{equ:asb}) in Appendix \ref{app:Lemma 5} using \eqref{equ:updae}. Then, there exist uniform bounds $L^k \leq \overline{L}, \theta_1^k \leq  \overline{\theta}_1^k, \theta_2^k\leq  \overline{\theta}_2, \theta_3^k\geq  \overline{\theta}_3$. This is because they are polynomials in $(\|A\|,\|B\|,\|Q\|,\|R\|, \underline{\sigma}(Q), \underline{\sigma}(R), \gamma^k, J_{\gamma^k}(K^k))$ that are uniformly bounded, i.e., $J_{\gamma^0}^*\leq J_{\gamma^k}(K^k)< 12(J_1^*)^2/((1-\xi)\underline{\sigma}(Q)), \gamma^0 \leq \gamma^k \leq 1$. 
	Then, $\tau, N_e, r, \eta$ in (\ref{condition}) are a sufficient condition for that in Lemma \ref{lem:grad} combining $D-J_{\gamma}^* > 2J_1^*-J_1^* = J_1^*$. To derive the expression of $M$, we provide an upper bound for the right-hand side of (\ref{equ:pgroll}) which at the $k$-th iteration is
	\begin{equation}\label{equ:iter_numb}
	-\log((J_{\gamma^k}(K^k)-J_{\gamma^k}^*)/(D - J_{\gamma^{k}}^*))/\log(1-\mu^k\eta^k/8),
	\end{equation}
	where $\mu^k  = \underline{\sigma}(R)/(2\|\Sigma^*(\gamma^k)\|)$ is the gradient dominance constant, and $\eta^k$ is the stepsize. Then, (\ref{equ:iter_numb}) can be upper bounded using the numerical inequality $\log(1+x)\leq x$ by
	\begin{equation}\label{equ:iter2}
	8\log((J_{\gamma^k}(K^k)-J_{\gamma^k}^*)/(D - J_{\gamma^{k}}^*))/(\mu^k\eta^k).
	\end{equation}
	Since $J_{\gamma^k}(K^k) \leq {12(J^*_1)^2}/{((1-\xi)\underline{\sigma}(Q))}$, the numerator of (\ref{equ:iter2}) can be upper bounded by 
	$$
	8\log\left(\frac{J_{\gamma^k}(K^k)-J_{\gamma^k}^*}{D - J_{\gamma^{k}}^*}\right) 
	< 8\log \left(\frac{12J^*_1}{(1-\xi)\underline{\sigma}(Q)}\right).
	$$
	For the denominator of (\ref{equ:iter2}), we note that the gradient dominance constant has a uniform lower bound $\mu^k \geq \underline{\sigma}(R)/(2J^*_1)$ for $k\in \mathbb{N}$, which follows from $J^*_1 \geq J_{\gamma^k}^* \geq \|\Sigma^*(\gamma^k)\|$. The stepsize $\eta^k $ is upper bounded by ${1}/{(\overline{w}\overline{L})}$ as indicated in (\ref{condition}). Then, the number of PG updates $M$ is uniformly bounded by 
	$$M \geq \frac{\overline{L} J^*_1\overline{\omega}}{\underline{\sigma}(R)} \log \left(\frac{12J^*_1}{(1-\xi)\underline{\sigma}(Q)}\right).$$
	By Lemma \ref{lem:grad}, the PG updates under the condition in (\ref{condition}) fail to return $J^k(K^{k+1}) \leq D$ at the $k$-th iteration with probability not larger than 
	$$
	c' M(n^{-\beta}+N_e^{-\beta}+N_ee^{-\frac{n}{8}}+e^{-c' N_e}).
	$$
	
	Finally, \eqref{equ:updae} holds by union bound for $k\in \mathbb{N}$ with probability in (\ref{equ:prob}).
\end{proof}

By Lemma \ref{lem:newtech}, the update rate $\alpha^k$ has a uniform lower bound  \eqref{equ:updae}. Thus, considering the effects of the coefficient $\xi$, the number of iterations for Algorithm \ref{alg:mfeda} to converge is 
$ {(6J_1^*-\underline{\sigma}(Q))}\log ({1}/{\gamma^0})/(\xi\underline{\sigma}(Q))$
with probability in (\ref{equ:prob}). The bound on the spectral radius of the resulted policy $\rho(A-BK^k)<\sqrt{1-(1-\xi)\underline{\sigma}(Q)/(6J_1^*)}$ follows from the derivation of (\ref{equ:spectralbound}) combining \eqref{equ:updae}.

\begin{IEEEbiography}[{\includegraphics[width=1in,height=1.25in,clip,keepaspectratio]{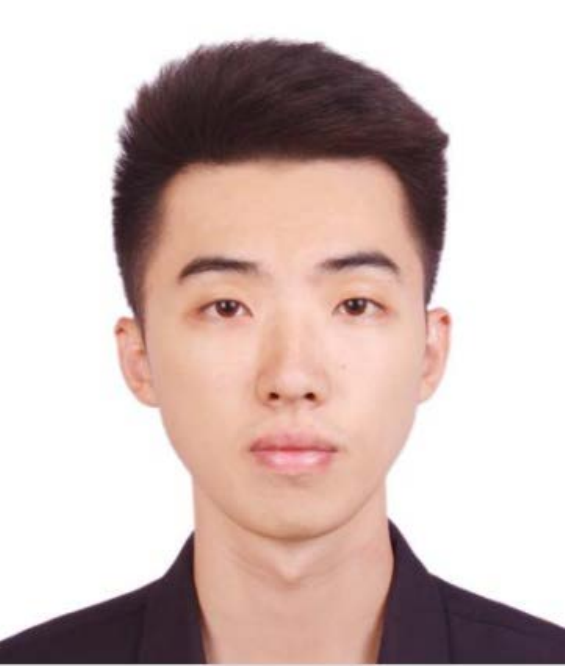}}]{Feiran Zhao} received the B.S. degree in Control Science and Engineering from the School of Astronautics, Harbin Institute of Technology, Harbin, China, in 2018. He is currently pursuing the Ph.D. degree in Control Science and Engineering at the Department of Automation, Tsinghua University, Beijing, China. His research interests include reinforcement learning, control theory and their intersections.
\end{IEEEbiography}

\begin{IEEEbiography}[{\includegraphics[width=1in,height=1.25in,clip,keepaspectratio]{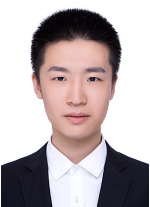}}]{Xingyun Fu} received the B.S. degree in Control Science and Engineering, in 2020, from the Department of Automation, Tsinghua University, Beijing, China, where he is currently pursuing the Ph.D. degree in Control Science and Engineering. His research interests include reinforcement learning and data-driven control.
\end{IEEEbiography}

\begin{IEEEbiography}
	[{\includegraphics[width=1in,height=1.25in,clip,keepaspectratio]{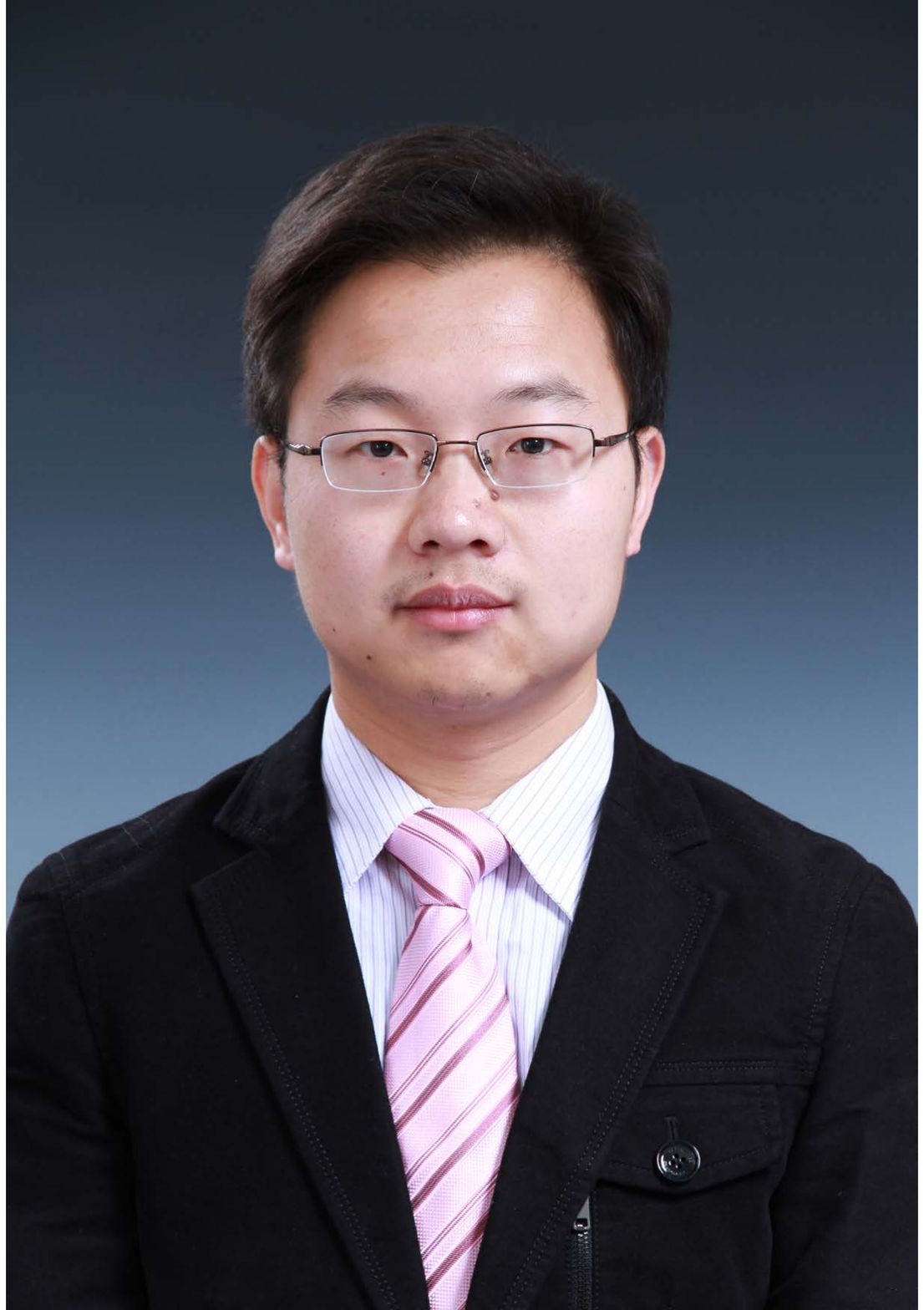}}]
	{Keyou You}(SM'17) received the B.S. degree in Statistical Science from Sun Yat-sen University, Guangzhou, China, in 2007 and the Ph.D. degree in Electrical and Electronic Engineering from Nanyang Technological University (NTU), Singapore, in 2012. After briefly working as a Research Fellow at NTU, he joined Tsinghua University in Beijing, China, where he is now a full professor in the Department of Automation. He held visiting positions at Politecnico di Torino,  Hong Kong University of Science and Technology,  University of Melbourne and etc. His  research interests include data-driven control, networked control systems, distributed optimization and learning, and their applications.
	
	Dr. You received the Guan Zhaozhi award at the 29th Chinese Control Conference in 2010 and the ACA (Asian Control Association) Temasek Young Educator Award in 2019. He has received the National Science Fund for Distinguished Young Scholars in 2023. He is currently an Associate Editor for Automatica,
	IEEE Transactions on Control of Network Systems, and IEEE Transactions on
	Cybernetics.
\end{IEEEbiography}
\end{document}